\documentclass{ws-ijmpc}

\usepackage{ulem}
\usepackage{hyperref}
\usepackage{lipsum,verbatim}
\usepackage{graphicx}
\usepackage{multirow,rotate,subfigure}
\usepackage{amscd,amsmath,amssymb}
\usepackage{booktabs}

\usepackage{adjustbox}

\usepackage{mathrsfs}
\usepackage{mathtools}

\usepackage{amscd,amssymb}
\usepackage{mathtools}

\usepackage{color, colortbl}
\definecolor{red}{rgb}{1,0,0}
\definecolor{green}{rgb}{0,1,0}
\definecolor{blue}{rgb}{0,0,1}

\definecolor{Gray}{gray}{0.9} 
\definecolor{LightCyan}{rgb}{0.88,1,1}
\definecolor{ciano}{rgb}{0,1,1}
\definecolor{magenta}{rgb}{1,0,1}
\definecolor{amarelo}{rgb}{1,1,0}
\definecolor{iw}{rgb}{0.7,0.93,0.36}  
\definecolor{lav}{rgb}{0.96,0.73,1.0}  
\newcommand{\iw}{\cellcolor{iw}}
\newcommand{\iy}{\cellcolor{lav}}

\begin{document}

\markboth{Jason A.C.~Gallas}
{Monogenic period equations are cyclotomic polynomials}

\catchline{}{}{}{}{}

\title{Monogenic period equations are cyclotomic polynomials}

\vspace{-0.25truecm}

\author{Jason A.C.~Gallas} 

\address{
Instituto de Altos Estudos da Para\'\i ba,
  Rua Silvino Lopes 419-2502,\\   
  58039-190 Jo\~ao Pessoa, Brazil,\\
Complexity Sciences Center, 9225 Collins Ave.~1208, 
Surfside FL 33154, USA,\\
  Max-Planck-Institut f\"ur Physik komplexer Systeme,
  01187 Dresden, Germany\\
  jgallas@pks.mpg.de}

\maketitle

\begin{history}
\received{19 December 2019}
\accepted{30 December 2019}
\centerline{https://doi.org/10.1142/S0129183120500588}
\end{history}


\begin{abstract}
We study monogeneity in  {\sl period equations}, $\psi_e(x)$,
the au\-xi\-li\-a\-ry equations introduced by
Gauss to solve cyclotomic polynomials by radicals.
All  monogenic $\psi_e(x)$ of degrees $4 \leq e \leq 250$ are
determined for extended intervals of primes $p=ef+1$,
and found to coincide either with cyclotomic polynomials,
or with simple de Moivre reduced forms of cyclotomic polynomials.
The former case occurs for $p=e+1$, and the latter for $p=2e+1$.
For $e\geq4$, we conjecture all monogenic period equations to
be cyclotomic polynomials.
Totally real period equations are of interest in
applications of quadratic discrete-time dynamical systems.
\end{abstract}

\ccode{PACS Nos.:
      02.70.Wz, 
      02.10.De, 
      03.65.Fd} 

\section{Introduction}

A recent paper in this Journal has shown that, in the partition generating
limit\cite{bk}, several orbital equations and clusters of orbital equations of the
quadratic (logistic) map coincide with cyclotomic period equations\cite{jg19}.
Period equations were introduced in 1801 by Gauss as {\it auxiliary equations}
to solve  cyclotomic polynomials by radicals.

The purpose of this paper is to report a startling finding obtained by
extensive empirical computations: we find {\sl monogenic} period equations to 
be either cyclotomic polynomials, or simple de Moivre reduced forms of cyclotomic
polynomials, thereby implying the existence of a hierarchical interdependence
among fields and subfields of cyclotomic polynomials, or orbital equations.
This fact is quite remarkable because, although cyclotomic polynomials are among
the most extensively studied polynomials for nearly 220 years, it seems to have
hitherto escaped attention that, in essence, Gauss auxiliary monogenic period
equations are nothing else than just cyclotomic polynomials themselves.

\section{Context and basic definitions}


As shown by Dedekind\cite{de78},
it is always possible to fix a number field $K$ of finite degree $n$ over $\mathbb Q$
by
selecting an algebraic integer $\alpha \in K$ such that $K=\mathbb Q(\alpha)$.
In other words, a number field $K$ may be determined
by selecting $\alpha$ as a root of a monic $n$-degree $\mathbb Q$-irreducible
polynomial $f(x)$ and expressing it in terms of $n$ integers
$\alpha_1, \alpha_2,\cdots,\alpha_n$, independent of each other, 
forming a {\sl basis}
for $\mathcal O_K$, the {\sl ring} of integers  in $K$.
In this context, a key problem
is to decide whether or not the ring  $\mathcal O_K$ is {\sl monogenic}, namely
if there exists an element $\alpha\in K$ such that $\mathcal O_K$ is a polynomial
ring $\mathbb Z[\alpha]$, i.e.~if  powers of the type
$(1, \alpha, \cdots, \alpha^{n-1})$ constitute a
{\sl power integral basis}\cite{g02,eg17,n04,ha64}.
Every algebraic number field has an integral basis, not necessarily a
power integral basis.

In a monogenic field $K$, the field discriminant $\Delta$ coincides with
the standard discriminant $D$ of the minimal polynomial of $\alpha$.
For non-monogenic fields such identity does not hold.
Generically,  $D$ and $\Delta$ are  
interconnected by the harmless-looking relation\cite{eg17,n04,ha64,kro}
\begin{equation}
  D=k^2 \Delta,   \label{discri}
\end{equation}  
for some $k\in\mathbb Z$, called {\sl index} by Dedekind\cite{de78}, and
``ausserwesentlicher Theiler der Discriminante'', {\sl inessential discriminant
divisor},  by Kronecker\cite{kro}.
As pointed out by Vaughan\cite{v85},
``{\sl while $D$ can be found by straightforward (if tedious) computation,
the value of $k$ is quite another story.
According to Cohn\cite{hc78}, page 77, for
example, to determine $k$, one would have to test a finite number (which
may be very large) of elements of $K$ to see if they are integral.}''
An additional complication to obtain $k$ is the fact that the choice
of $\alpha$ is not unique and, therefore, there are several distinct minimal
polynomials from which to compute $D$.
So, one may consider $k$ as a sort of ``quality measure''
for minimal polynomial representation and for monogeneity.
As described in \S 4 of Dedekind's paper\cite{de78}, for quite some time he
believed to be always possible to find a suitable $\alpha$ leading to a power basis.
This, 
until he found what became a popular textbook example of a non-monogenic field
generated
by a root of $x^3 -x^2 -2x-8=0$, for which $\Delta=-503$, $k^2=4$, and $D=-4\cdot503$. 

The computation of the index
$k$ and the verification whether or not a given number field
has a power basis are two hard problems\cite{g02,eg17}.
A taste for the difficulties and the representative computation times involved
in such tasks may be obtained from a paper by Bilu et al.~\cite{bi04}.
Efficient
algorithms for determining generators of power integral bases involve solving
Diophantine equations known as {\sl index form equations}\cite{eg17}.

The first general algorithm for determining all power integral bases in
number fields was given in 1976 by Gy\H{o}ry\cite{k76}.
Subsequently, efficient algorithms were elaborated to determine
power integral bases for number fields of degree at most six and some
special classes of higher degree number fields\cite{g02,eg17,bi04}.
Section 7.3 of the book by Evertse and Gy\H{o}ry\cite{eg17}
discusses significant results by Gras on Abelian number fields
of degree $n$, where $n$ is relatively prime to 6. See also Ref.~\cite{sh14}.
Finally, we mention a general and still largely open problem stated by Hasse:
to give a characterization of monogenic number fields\cite{eg17,n04,ha64}.
Hasse's problem has been considered, among others, by Nakahara and
co-authors in several contexts during the last 50 years or so.
See Refs.~\cite{sh14,hn15,ka16,ka19,sn15} and references therein.
See also Evertse\cite{eee}.

This paper reports results of an extended investigation of the distribution
of monogeneity in cyclotomic period equations $\psi_e(x)$, a wide class of
functions  underlying the solution of cyclotomic polynomials\cite{da,bach,reu}.
{ We find that all monogenic period equations are either cyclotomic polynomials
or simple reduced forms of cyclotomic polynomials.}
Here, monogenic period equations are obtained with the help of an
expression for the field discriminant $\Delta_e$
of period equations\cite{jg19},  Eq.~(\ref{delta}) below.
Since the polynomial discriminant $D$ may be easily computed,
knowledge of $\Delta_e$ gives at once  
\begin{equation}
    k^2 = {D}/{\Delta_e},  \label{kkk}
\end{equation}
which is a convenient tool to sort out all equations with index $k=1$.
In what follows, we present results obtained for such monogenic period equations.
Equations (\ref{kkk}) and (\ref{delta}) grant access to families of
equations of arbitrarily high degrees $e$, opening the possibility of
studying monogeneity well beyond the aforementioned low-degree limits. 
Note that for high degrees the division in Eq.~(\ref{kkk})
involves exceedingly large integers.

\section{Field discriminants of cyclotomic period equations}
\label{pereq}

Let $g$ be a primitive root of a prime $p=ef+1$, and 
$r=\exp(2\pi i/p)$. In the {\sl Disquisitiones Arithmetic\ae},
Gauss defined $e$  sums $\eta_i$
called ``periods''\cite{da,bach,reu,jg19,ra64}:
\begin{equation*}
  \eta_i = \sum_{k=0}^{f-1} r^{g^{ke+i}}, \qquad i=0,1,\cdots,e-1,
\end{equation*}
With them, he defined {\sl period equations} $\psi_e(x)$,
polynomials of degree $e$ whose roots are the periods $\eta_i$  
\begin{equation*}
  \psi_e(x) = \prod_{k=0}^{e-1} (x-\eta_{k})
            = x^e+x^{e-1}+\alpha_2x^{e-2}+\cdots+ \alpha_e,
               \qquad \alpha_\ell \in \mathbb Z. \label{theta}
\end{equation*}

Period equations $\psi_e(x)$ constitute a wide class of equations for
which the computation of the field discriminant $\Delta_e$
presents no difficulties, being given by\cite{jg19}
\begin{equation}\label{delta}
    \Delta_e =
    \begin{cases*}
  -p^{e-1}, & {\rm if} $(e-1) \hbox{ \rm mod } 4 =1 {\ \ \ \rm and\ \ \ }
                       f {\ \rm mod\ } 2 = 1 $, \\
  \phantom{-}p^{e-1}, & {\rm if otherwise}.
    \end{cases*}
\end{equation}
Together with Eq.~(\ref{kkk}), this discriminant provides a handy
criterion to sort out $k^2=1$ monogenic equations through a simple
division of two (possibly very large) integers.

\section{Properties of monogenic period equations}

Table \ref{tab:tab01}
lists monogenic period equations as a
function of $e$ for the first few equations of a much longer list,
containing seven equations for every $e\leq 250$.
The table also displays the signature of $\psi_e(x)$.
The {\sl signature}\cite{ha64} of a polynomial is
the doublet $(n_R, n_P)$, sometimes written more economically
as $n_R$, 
informing the number $n_R$ of real roots of $\psi_e(x)$, and the
number $n_P$ of {\it pairs} of complex roots. 
Table  \ref{tab:tab01} illustrates regularities that are
consistently observed for $e$ up to 250.

Among the equations obtained for a given $e$ we find no more
than two types of polynomials leading to  $k=1$. They have
either  totally complex roots, $n_R=0$, or totally real, $n_R=e$.
Polynomials with  $n_R=0$ are highlighted differently,
to reflect the sign of their discriminants.

For $e=2$, period equations are quadratic and, as known, are all
monogenic\cite{g02}.
For $e=3$, we determined the  growth of the number of $k=1$ equations
as a function of $e$, up to $e=9000$. Such growth obeys
a power-law distribution, implying  the existence of
an infinite number of monogenic cubic equations.

For $e\geq4$, we find no more than two monogenic equations for each
value of $e$, as illustrated in Table \ref{tab:tab01}.
From the Table one recognizes a trend observed also for higher
values of $e$: the absence of monogenic equations for several
values of $e$.
For instance, for  $e\leq 100$ we find no monogenic period equations
for
$e=7$, $13$, $17$, $19$, $24$, $25$, $27$, $31$, $32$, $34$,
$37$, $38$, $43$, $45$, $47$, $49$, $55$, $57$, $59$, $61$,
$62$, $64$, $67$, $71$, $73$, $76$, $77$, $79$, $80$, $84$, $85$,
$87$, $91$, $92$, $93$, $94$, and $97$.
Analogously, there are 62 cases of missing cyclotomic polynomials
with degree $\leq 100$.

Period equations are not difficult to generate fast and explicitly
up to very high degrees. As already mentioned,
this means that Eqs.~(\ref{kkk}) and (\ref{delta}) open the possibility
to investigate monogeneity systematically for an important family of
equations well beyond the aforementioned low-degree limits.

\begin{table}[!tb]
\tbl{Monogenic period equations  as a function of $e$
  for primes $p=ef+1$ and signature $n_R$.
Highlighted equations have signature $n_R=0$ and coincide
with cyclotomic polynomials $\Phi_p(x)$.
Non-highlighted equations and signature $n_R=e$ are
{\it reduced\/} cyclotomic polynomials.}  
{\begin{tabular}{@{}|c|c|c|c|l|@{}} \toprule        
\hline      
  $e$ & $p$ & $n_R$ & $D=\Delta_e$  & $\psi_e(x)$ \\
\hline
\iw4& \iw5 & \iw0&\iw $5^3$ & \iw ${x}^{4}+{x}^{3}+{x}^{2}+x+1$\\
5& 11 & 5&  $11^4$ & ${x}^{5}+{x}^{4}-4\,{x}^{3}-3\,{x}^{2}+3\,x+1$\\
\iy6 & \iy7 &\iy 0&\iy  $-7^5$  & \iy ${x}^{6}+{x}^{5}+{x}^{4}+{x}^{3}+{x}^{2}+x+1$\\
6& 13 & 6&  $13^5$ & ${x}^{6}+{x}^{5}-5\,{x}^{4}-4\,{x}^{3}+6\,{x}^{2}+3\,x-1$\\
8& 17 & 8& $17^7$  & ${x}^{8}+{x}^{7}-7\,{x}^{6}-6\,{x}^{5}+15\,{x}^{4}+10\,{x}^{3}-10\,{x}^{2}-4\,x+1$\\
9& 19 & 9& $19^8$ & ${x}^{9}+{x}^{8}-8\,{x}^{7}-7\,{x}^{6}+21\,{x}^{5}+15\,{x}^{4}-20\,{x}^{3}-10\,{x}^{2}+5\,x+1$\\
\iy10&\iy 11&\iy 0&\iy $-11^9$ & \iy ${x}^{10}+{x}^{9}+{x}^{8}+{x}^{7}+{x}^{6}+{x}^{5}+{x}^{4}
                               +{x}^{3}+{x}^{2}+x+1$\\
11& 23& 11& $23^{10}$&${x}^{11}+{x}^{10}-10\,{x}^{9}-9\,{x}^{8}+36\,{x}^{7}+28\,{x}^{6}-56\,{x}^{5}-35\,{x}^{4}+35\,{x}^{3}+15\,{x}^{2}-6\,x-1$\\
\iw12&\iw 13&\iw 0&\iw $13^{11}$& \iw ${x}^{12}+{x}^{11}+{x}^{10}+{x}^{9}+{x}^{8}+{x}^{7}+{x}^{6}+{x}^{5}+{x}^{4}+{x}^{3}+{x}^{2}+x+1$\\
14& 29& 14& $29^{13}$ & ${x}^{14}+{x}^{13}-13\,{x}^{12}-12\,{x}^{11}+66\,{x}^{10}+55\,{x}^{9}-165\,{x}^{8}-120\,{x}^{7}$\cr
 &&&& \qquad $+210\,{x}^{6}+126\,{x}^{5}-126\,{x}^{4}-56\,{x}^{3}+28\,{x}^{2}+7
\,x-1$\\
15& 31& 15& $31^{14}$ & ${x}^{15}+{x}^{14}-14\,{x}^{13}-13\,{x}^{12}+78\,{x}^{11}+66\,{x}^{10}-220\,{x}^{9}-165\,{x}^{8}+330\,{x}^{7}$\cr
&&&& \qquad $+210\,{x}^{6}-252\,{x}^{5}-126\,{x}^{4}+84\,{x}^{3}
+28\,{x}^{2}-8\,x-1$\\
\iw16&\iw 17&\iw 0&\iw $17^{15}$ & \iw ${x}^{16}+{x}^{15}+{x}^{14}+{x}^{13}+{x}^{12}+{x}^{11}+{x}^{10}+{x}^{9}+{x}^{8}+{x}^{7}$\cr
&&&& \qquad \iw$+{x}^{6}+{x}^{5}+{x}^{4}+{x}^{3}+{x}^{2}+x+1$\\
\iy18&\iy 19&\iy 0&\iy $-19^{17}$ & \iy ${x}^{18}+{x}^{17}+{x}^{16}+{x}^{15}+{x}^{14}+{x}^{13}+{x}^{12}+{x}^{11}+{x}^{10}+{x}^{9}+{x}^{8}+{x}^{7}$\cr
&&&&\qquad \iy $+{x}^{6}+{x}^{5}+{x}^{4}+{x}^{3}+{x}^{2}+x+1$\\
18& 37& 18& $37^{17}$ & ${x}^{18}+{x}^{17}-17\,{x}^{16}-16\,{x}^{15}+120\,{x}^{14}+105\,{x}^{13}-455\,{x}^{12}-364\,{x}^{11}+1001\,{x}^{10}$\cr
&&&&\qquad $+715\,{x}^{9}-1287\,{x}^{8}-792\,{x}^{7}+924
\,{x}^{6}+462\,{x}^{5}-330\,{x}^{4}-120\,{x}^{3}+45\,{x}^{2}+9\,x-1$\\
20& 41& 20& $41^{19}$ & ${x}^{20}+{x}^{19}-19\,{x}^{18}-18\,{x}^{17}+153\,{x}^{16}+136\,{x}^{15}-680\,{x}^{14}-560\,{x}^{13}+1820\,{x}^{12}$\cr
&&&&\qquad $+1365\,{x}^{11}-3003\,{x}^{10}-2002\,{x}^{9}+
3003\,{x}^{8}+1716\,{x}^{7}-1716\,{x}^{6}-792\,{x}^{5}$\cr
&&&&\qquad $+495\,{x}^{4}+165\,{x}^{3}-55\,{x}^{2}-10\,x+1$\\
21& 43& 21& $43^{20}$ & ${x}^{21}+{x}^{20}-20\,{x}^{19}-19\,{x}^{18}+171\,{x}^{17}+153\,{x}^{16}-816\,{x}^{15}-680\,{x}^{14}+2380\,{x}^{13}$\cr
&&&&\qquad $+1820\,{x}^{12}-4368\,{x}^{11}-3003\,{x}^{10}
+5005\,{x}^{9}+3003\,{x}^{8}-3432\,{x}^{7}-1716\,{x}^{6}$\cr
&&&&\qquad $+1287\,{x}^{5}+495\,{x}^{4}-220\,{x}^{3}-55\,{x}^{2}+11\,x+1$\\
\iy22&\iy 23&\iy 0&\iy $-23^{21}$ &\iy ${x}^{22}+{x}^{21}+{x}^{20}+{x}^{19}+{x}^{18}+{x}^{17}+{x}^{16}+{x}^{15}+{x}^{14}+{x}^{13}+{x}^{12}+{x}^{11}+{x}^{10}$\cr
&&&&\qquad \iy $+{x}^{9}+{x}^{8}+{x}^{7}+{x}^{6}+{x}^{5}+{x
}^{4}+{x}^{3}+{x}^{2}+x+1$\\
23& 47& 23& $47^{22}$ & ${x}^{23}+{x}^{22}-22\,{x}^{21}-21\,{x}^{20}+210\,{x}^{19}+190\,{x}^{18}-1140\,{x}^{17}-969\,{x}^{16}+3876\,{x}^{15}$\cr
&&&&\qquad $+3060\,{x}^{14}-8568\,{x}^{13}-6188\,{x}^{12
}+12376\,{x}^{11}+8008\,{x}^{10}-11440\,{x}^{9}-6435\,{x}^{8}$\cr
&&&&\qquad $+6435\,{x}^{7}+3003\,{x}^{6}-2002\,{x}^{5}-715\,{x}^{4}+286\,{x}^{3}+66\,{x}^{2}-12\,x-1$\\
26& 53& 26& $53^{25}$ & ${x}^{26}+{x}^{25}-25\,{x}^{24}-24\,{x}^{23}+276\,{x}^{22}+253\,{x}^{21}-1771\,{x}^{20}-1540\,{x}^{19}+7315\,{x}^{18}$\cr
&&&&\qquad $+5985\,{x}^{17}-20349\,{x}^{16}-15504\,{x}^
{15}+38760\,{x}^{14}+27132\,{x}^{13}-50388\,{x}^{12}$\cr
&&&&\qquad $-31824\,{x}^{11}+43758\,{x}^{10}+24310\,{x}^{9}-24310\,{x}^{8}-11440\,{x}^{7}+8008\,{x}^{6}$\cr
&&&&\qquad $+3003\,{x}^{5}-1365\,{x}^{4}-364\,{x}^{3}+91\,{x}^{2}+13\,x-1$\\
\iw28 &\iw 29&\iw 0&\iw $29^{27}$ &\iw ${x}^{28}+{x}^{27}+{x}^{26}+{x}^{25}+{x}^{24}+{x}^{23}+{x}^{22}+{x}^{21}+{x}^{20}+{x}^{19}+{x}^{18}+{x}^{17}+{x}^{16}$\cr
&&&&\qquad\iw $+{x}^{15}+{x}^{14}+{x}^{13}+{x}^{12}+{x}^{11}+{x}^{10}+{x}^{9}+{x}^{8}+{x}^{7}+{x}^{6}+{x}^{5}+{x}^{4}+{x}^{3}+{x}^{2}+x+1$\\
29& 59& 29& $59^{28}$ & ${x}^{29}+{x}^{28}-28\,{x}^{27}-27\,{x}^{26}+351\,{x}^{25}+325\,{x}^{24}-2600\,{x}^{23}-2300\,{x}^{22}+12650\,{x}^{21}$\cr
&&&&\qquad $+10626\,{x}^{20}-42504\,{x}^{19}-33649\,{x}^{18}+100947\,{x}^{17}+74613\,{x}^{16}-170544\,{x}^{15}$\cr
&&&&\qquad $-116280\,{x}^{14}+203490\,{x}^{13}+125970\,{x}^{12}-167960\,{x}^{11}-92378\,{x}^{10}+92378\,{x}^{9}$\cr
&&&&\qquad $+43758\,{x}^{8}-31824\,{x}^{7}-12376\,{x}^{6}+6188\,{x}^{5}+1820\,{x}^{4}-560\,{x}^{3}-105\,{x}^{2}+15\,x+1$\\
\iy30&\iy 31&\iy 0&\iy $-31^{29}$ &\iy ${x}^{30}+{x}^{29}+{x}^{28}+{x}^{27}+{x}^{26}+{x}^{25}+{x}^{24}+{x}^{23}+{x}^{22}+{x}^{21}+{x}^{20}+{x}^{19}+{x}^{18}+{x}^{17}+{x}^{16}$\cr
&&&&\qquad\iy $+{x}^{15}+{x}^{14}+{x}^{13}+{x}^{12}+{x}^{11}+{x}^{10}+{x}^{9}+{x}^{8}+{x}^{7}+{x}^{6}+{x}^{5}+{x}^{4}+{x}^{3}+{x}^{2}+x+1$\\
30& 61& 30& $61^{29}$ & ${x}^{30}+{x}^{29}-29\,{x}^{28}-28\,{x}^{27}+378\,{x}^{26}+351\,{x}^{25}-2925\,{x}^{24}-2600\,{x}^{23}+14950\,{x}^{22}$\cr
&&&&\qquad $+12650\,{x}^{21}-53130\,{x}^{20}-42504\,{x}^{19}+134596\,{x}^{18}+100947\,{x}^{17}-245157\,{x}^{16}$\cr
&&&&\qquad $-170544\,{x}^{15}+319770\,{x}^{14}+203490\,{x}^{13}-293930\,{x}^{12}-167960\,{x}^{11}+184756\,{x}^{10}$\cr
&&&&\qquad $+92378\,{x}^{9}-75582\,{x}^{8}-31824\,{x}^{7}+18564\,{x}^{6}+6188\,{x}^{5}-2380\,{x}^{4}-560\,{x}^{3}+120\,{x}^{2}+15\,x-1$\\
33& 67& 33& $67^{32}$ & ${x}^{33}+{x}^{32}-32\,{x}^{31}-31\,{x}^{30}+465\,{x}^{29}+435\,{x}^{28}-4060\,{x}^{27}-3654\,{x}^{26}+23751\,{x}^{25}+20475\,{x}^{24}$\cr
&&&&\qquad $-98280\,{x}^{23}-80730\,{x}^{22}+296010\,{x}^{21}+230230\,{x}^{20}-657800\,{x}^{19}-480700\,{x}^{18}+1081575\,{x}^{17}$\cr
&&&&\qquad $+735471\,{x}^{16}-1307504\,{x}^{15}-817190\,{x}^{14}+1144066\,{x}^{13}+646646\,{x}^{12}-705432\,{x}^{11}$\cr
&&&&\qquad $-352716\,{x}^{10}+293930\,{x}^{9}+125970\,{x}^{8}-77520\,{x}^{7}-27132\,{x}^{6}+11628\,{x}^{5}+3060\,{x}^{4}$\cr
&&&&\qquad $-816\,{x}^{3}-136\,{x}^{2}+17\,x+1$\\
\hline
\end{tabular}}\label{tab:tab01}
\end{table}

\section{Monogenic period equations are cyclotomic polynomials}

Gauss showed how to decompose and to solve explicitly 
in terms of radicals cyclotomic polynomials $x^p-1=0$, for prime $p$.
The procedure is described in \S 343 of his
{\sl Disquisitiones Arithmetic\ae}\cite{da,bach,reu,ra64}.
To this end, he introduced period equations as mere
{\sl auxiliary equations}.
This fact notwithstanding, by examining Table \ref{tab:tab01} we
deduce that
{\sl   
  all monogenic period equations $\psi_e(x)$ are
  nothing else than cyclotomic polynomials $\Phi_p(x)$,
  interconnected either by $p=e+1$ or $p=2e+1$}.
This latter interconnection may be identified as follows.

Consider, e.g., $e=5$  and $p=2\cdot5+1=11$, when the corresponding
cyclotomic polynomial is
\[ \Phi_{11}(x)=\frac{x^{11}-1}{x-1} = x^{10} + x^{9} + \cdots + x^2+x+1. \]
By changing variables according to a standard
de Moivre transformation\cite{dm30}
\begin{equation}
  z=x+1/x,  \label{direta}
\end{equation}  
and replacing $z$ by $x$ in the equation so obtained, one finds
\[ \psi_5(x) = {x}^{5}+{x}^{4}-4\,{x}^{3}-3\,{x}^{2}+3\,x+1. \]
This quintic, listed in Table \ref{tab:tab01},
was solved explicitly by radicals in a {\sl M\'emoire}
read in November 1770 by Vandermonde\cite{v70,i88}.

Conversely, changing $x$ in $\psi_5(x)$ according to the
{\sl dual transformation}
\begin{equation}
  x=z+1/z  \label{inversa}
\end{equation}
(and replacing $z$ by $x$) one recovers the cyclotomic $\Phi_{11}(x)$.
Thus,  it is an easy matter to pass from one polynomial to the other one,
showing that, essentially, monogenic period equations are 
cyclotomic polynomials. Here is another example.

In 1796, almost three decades after Vandermonde\cite{l18},
Gauss recorded in his
mathematical diary\cite{k03,Tbuch,du04} that the regular 17-gon
can be constructed by ruler and compass alone.
In print, his construction appeared\cite{da} only in 1801.
The solution amounts to reducing by two, four times in succession,
the degree of the  17-gon cyclotomic polynomial, namely
\[ \Phi_{17}(x)=\frac{x^{17}-1}{x-1} = x^{16} + x^{15} + \cdots + x^2+x+1. \]
Applying de Moivre's transformation, Eq.~(\ref{direta}), to $\Phi_{17}(x)$
one gets the first of such reductions,
also present in Table \ref{tab:tab01}:\quad
\[\psi_8(x) = {x}^{8}+{x}^{7}-7\,{x}^{6}-6\,{x}^{5}+15\,{x}^{4}
               +10\,{x}^{3}-10\,{x}^{2}-4\,x+1. \]

Many additional examples are obtained in
a similar way, by using the dual transformation, Eq.~(\ref{inversa}),
to unfold $\psi_e(x)$ with signature $n_R=e$
for all equations listed in Table \ref{tab:tab01},
thereby obtaining the associated cyclotomic $\Phi_{2e+1}(x)$.
The dual transformation worked also for all additional equations up
to $e=250$ (not all in Table \ref{tab:tab01}).
For instance,  $e=96$ is the largest value $\leq100$ with two
monogenic period equations. For signature $n_R=0$ and $f=1$ we find
\begin{alignat*}{2}
  \psi_{96}(x) &= x^{96}+x^{95}+x^{94}+x^{93}+x^{92}+\cdots+x^4+x^3+x^2+x+1,
\end{alignat*}
while for  $n_R=96$ and $f=2$ we find
\begin{alignat*}{2}
  \psi_{96}(x) &= x^{96}+x^{95}-95x^{94}-94x^{93}+4371x^{92}+4278x^{91}-129766x^{90}+\cdots\cr
   &\quad  -18009460x^6-2118760x^5+230300x^4+18424x^3-1176x^2-48x+1.
\end{alignat*}    
Similar doublets occur for
$e=6, 18, 30, 36, 78, 96, 138, 156, 198, 210, 228$, $270$, $306$, $330, \dots$.
There are 187 doublets for
$e\leq10^4$, 1164 for $e\leq10^5$, 7750 for $e\leq10^6$, etc.

\section{Conclusions}

The compelling computational evidence reported here leads us to conjecture
that for $e\geq4$ there are two classes of coincidences between monogenic
period equations $\psi_e(x)$ and cyclotomic polynomials $\Phi_p(x)$
interconnected by $p=ef+1$: \ 
The class of totally complex period equations, for which $p=e+1$,
and the class of totally real period equations, for which $p=2e+1$.
For all other values of $f$ in these classes, we only found non-monogenic
period equations and no connections to cyclotomic polynomials.
For $e=3$, as already mentioned, it is possible to find  an apparently
unbounded supply of monogenic period equations with $f>2$.
Totally real period equations are of significant interest for applications
in quadratic discrete-time dynamical systems in the partition
generating limit\cite{bk,jg19,g19,eg06,eg06b}.

\section*{Acknowledgments}
The author would like to thank
 Prof.~K.~Gy\H{o}ry, Debrecen, and Prof.~W.~Narkiewicz,
Wroc\l aw, for their kind feedback and helpful suggestions.
This work was supported by the Max-Planck Institute for the Physics of
Complex Systems, Dresden, in the framework of the Advanced Study Group on
{\sl Forecasting with Lyapunov vectors}.
The author was supported by CNPq, Brazil.





\end{document}